\documentclass[12pt]{article}
\usepackage{amscd, amssymb, latexsym, amsmath}

\newtheorem{thm}{Theorem}

\newtheorem{pro}{Proposition}

\newtheorem{dfn}{Definition}

\newtheorem{rem}{Remark}
\numberwithin{thm}{section}
\numberwithin{cor}{section}
\numberwithin{pro}{section}
\numberwithin{dfn}{section}
\numberwithin{lem}{section}
\numberwithin{rem}{section}
\numberwithin{equation}{section}

\newcommand{\f}{\frac}

\def\sm{\setminus}

\newcommand{\ba}{\begin{array}}
\newcommand{\ea}{\end{array}}

\newcommand{\ra}{\rightarrow}
\newcommand{\R}{\mathbb R}
\newcommand{\C}{\mathbb C}

\newcommand{\dvt}{d|\!|V_t|\!|}
\newcommand{\dvzero}{d|\!|V_0|\!|}

\def\eq#1{{\rm(\ref{#1})}}

\def\ts{\textstyle}

\def\Im{\mathop{\rm Im}}
\def\la{\lambda}
\def\th{\theta}
\def\al{\alpha}
\def\ov{\overline}

\def\de{\delta}
\def\ep{\epsilon}
\def\ga{\gamma}
\begin{document}
\title{Hamiltonian stationary cones and self-similar solutions in higher
dimension}
 \author{Yng-Ing Lee* and Mu-Tao Wang**}

\date{June 07, 2007, last revised February 3, 2008}
\maketitle \leftline{*Department of Mathematics and Taida
Institute of Mathematical Sciences,}
\leftline{\text{ }\,National
Taiwan University, Taipei, Taiwan}
 \leftline{\text{ }\,National
Center for Theoretical Sciences, Taipei Office}
 \centerline{email: yilee@math.ntu.edu.tw}
\leftline{**Department of Mathematics, Columbia University,  New
York, NY 10027, USA} \centerline{email: mtwang@math.columbia.edu}

\begin{abstract}

In \cite{lw}, we construct examples of two-dimensional Hamiltonian
stationary self-shrinkers and self-expanders for Lagrangian mean
curvature flows, which are asymptotic to the union of two
Schoen-Wolfson cones.  These self-shrinkers and self-expanders can
be glued together to yield solutions of the Brakke flow - a weak
formulation of the mean curvature flow.  Moreover, there  is no
mass loss along the Brakke flow. In this paper, we generalize
these results to higher dimension. We  construct new higher
dimensional Hamiltonian stationary cones of different topology as
generalizations of the Schoen-Wolfson cones. Hamiltonian
stationary self-shrinkers and self-expanders that are asymptotic
to these Hamiltonian stationary cones are also constructed.  They
can also be glued together to produce eternal solutions of the
Brakke flow without mass loss. Finally, we show the same conclusion holds for those Lagrangian
self-similar examples recently found by Joyce, Tsui and the first author
in \cite{jlt}.
\end{abstract}
\section{Introduction}

The existence of special Lagrangians in  Calabi-Yau manifolds
received much attention recently due to the critical role it plays
in the T-duality formulation of Mirror symmetry of
Strominger-Yau-Zaslow \cite{syz}. Special Lagrangians are
calibrated submanifolds and thus are volume minimizers \cite{hl}.
One potential approach to the construction of special Lagrangians
is the mean curvature flow- as the negative gradient flow of the
volume functional. However, the long-time existence of such flows
can only be verified in some special cases, see for example
\cite{sm}, \cite{swa}, \cite{wa1}, and \cite{wa2}. In this
article, we construct special weak solutions of the Lagrangian
mean curvature flows.

Our ambient space is always the complex Euclidean space $\C^n$
with coordinates $z^i=x^i+\sqrt{-1} y^i$, the standard symplectic
form $\omega= \sum_{i=1}^n dx^i \wedge dy^i$, and the standard
almost complex structure $J$ with $J(\frac{\partial}{\partial
x^i})=\frac{\partial}{\partial y^i}$. A Lagrangian submanifold is
an n-dimensional submanifold in $\C^n$, on which the symplectic
form $\omega$ vanishes. On a Lagrangian submanifold $L$, the mean
curvature vector ${H}$ is given by
\begin{equation}\label{mcv} {H}=J\nabla \theta\end{equation} where
$\theta$ is the Lagrangian angle and $\nabla$ is the gradient on
$L$. The Lagrangian angle $\theta$ can be defined by the relation
that
\[*_L (dz^1\wedge\cdots\wedge dz^n)=e^{i\theta}\] where
$*_L$ is the Hodge *-star operator on $L$. We recall

\begin{dfn}
A Lagrangian submanifold $L$ is called Hamiltonian stationary if
the Lagrangian angle is harmonic. i.e. $\Delta \theta=0$ where
$\Delta$ is the Laplace operator on $L$. $L$ is a special
Lagrangian if $\theta$ is a constant function.
\end{dfn}

A Hamiltonian stationary Lagrangian submanifold is a critical
point of the volume functional among all Hamiltonian deformations
and a special Lagrangian is a volume minimizer in its homology
class.

 By the first variation formula, the mean curvature vector
points to the direction where the volume is decreased most
rapidly. As the special Lagrangians are volume minimizers, it is
thus natural
 to use the mean curvature flow in the construction of special
 Lagrangians. Equation (\ref{mcv}) implies
that the mean curvature flow is a Lagrangian deformation, i.e. a
Lagrangian submanifold remains Lagrangian along the mean curvature
flow. In a geometric flow, the singularity often models on a
soliton solution. In the case of mean curvature flows , one type
of soliton solutions of particular interest are those moved by
scaling in the Euclidean space. We recall:

\begin{dfn}
A submanifold of the Euclidean space is called a
self-similar solution if

\[F^\perp=2cH\] for some nonzero constant $c$, where $F^\perp$ is normal
projection of the position vector $F$ in the Euclidean space and
$H$ is the mean curvature vector. It is called a self-shrinker
if $c<0$  and self-expander if $c>0$.
\end{dfn}

It is not hard to see that if $F$ is a self-similar solution, then
$F_t$ defined by $F_t=\sqrt{\frac{t}{c}} F$ is moved by the mean
curvature flow. By Huisken's monotonicity formula \cite{hu}, any
central blow up of a finite-time singularity of the mean curvature
flow is a self-similar solution.  In this article, we obtain
higher dimensional Hamiltonian stationary cones with different
topology as generalizations of the Schoen-Wolfson cones. We also
obtain Hamiltonian stationary self-shrinkers and self-expanders
which are asymptotic to these cones. Altogether they form
solutions of the Brakke flow (see \S \ref{brakke} ) which is a
weak formulation of the mean curvature flow proposed by Brakke in
\cite{br}.   To be more precise, we prove:
\begin{thm}\label{mainthm}
Assume $\lambda_j$ are non-zero integers and
$\sum_{j=1}^n\lambda_j>0$. Define
\[ \begin{split} V_t=&\{(x_1 e^{i\lambda_1s}, \cdots , x_n
e^{i\lambda_ns}):\quad 0\leq s < \pi,\\&\qquad \sum_{j=1}^n
\lambda_jx_j^2=(-2t)\sum_{j=1}^n \lambda_j, (x_1, \cdots, x_n)\in
\R^n \}.\end{split}\] Then $V_t$ is Hamiltonian stationary. It is
a self-shrinker for $t<0$, a  self-expander for $t>0$, and a cone
for $t=0$. Moreover, the varifold $\cup_t V_t$, $-\infty <t<
\infty$, form an eternal solution for Brakke flow without mass
loss.
\end{thm}
Our construction of Lagrangian self-similar solutions is
generalized to the non-Hamiltonian stationary case by Joyce, Tsui
and the first author in \cite{jlt}. These examples can also be
glued together to yield eternal solutions of the Brakke flow
without mass loss.

\begin{thm}\label{thm2} Let $k$ be a positive integer less than $n$.
Given $\lambda_j>0$ for $1\leq j\leq k <n$ and $\lambda_j<0$ for
$k< j\le n,$ let $w_1(s),\cdots,w_n(s):\R\ra\C\sm\{0\}$ be those
periodic functions with period $T$, which are obtained in {\rm
Theorem F} in {\rm\cite{jlt}} with $\al=1$. Define
\[ V_t=\{(x_1 w_1(s), \cdots , x_n
w_n(s)):\, 0\leq s < T,\, \sum_{j=1}^n \lambda_jx_j^2=2t, (x_1,
\cdots, x_n)\in \R^n \}.\] Then $V_t$ is  a  Lagrangian
self-shrinker for $t<0$, a Lagrangian self-expander for $t>0$, and
a Lagrangian cone for $t=0$. Moreover, the varifold $\cup_t V_t$,
$-\infty <t< \infty$, form an eternal solution for Brakke flow
without mass loss.
\end{thm}

The choice of $\al=1$ in Theorem~\ref{thm2} is arbitrary; there are a lot of freedom to rescale the
constants as discussed in \cite[Remark 3.2]{jlt}. Indeed, choosing $\al=-1$
instead will perhaps makes the statement  more
consistent with Theorem~\ref{mainthm}.

Theorem~\ref{mainthm} and \ref{thm2} are analogous to the
Feldman-Ilmanen-Knopf \cite{fik} gluing construction for the
K\"ahler-Ricci flows. Unlike the mean curvature flow, a notion of
weak solutions of Ricci flow has not yet been established.

This article is organized as the follows. The Hamiltonian
stationary examples in Theorem \ref{mainthm} and their geometry and topology are presented
in \S2. In \S3 , we recall the formulation of Brakke
flow  and prove Theorem \ref{mainthm}. The proof of
Theorem~\ref{thm2} and the discussion of the geometric properties of these examples are in \S4.

The first author would like to thank R. Schoen  for helpful
discussions and  hospitality during her visit in Stanford
University. The second author wishes to thank the support of the
Taida Institute for Mathematical Sciences during the preparation
of this article. The first author is supported by Taiwan NSC grant
96-2628-M-002. The second author is supported by NSF grant
DMS0605115 and a Sloan research fellowship.

\section{Hamiltonian Stationary Examples}
\subsection{The constructions}
For any $n$ nonzero integers $\lambda_1,\cdots, \lambda_n$,
consider the submanifold $L$ of $\C^n$ defined by

\[\{(x_1 e^{i\lambda_1s}, \cdots , x_n e^{i\lambda_ns})\,|
\,\, 0\leq s < 2\pi,\,\sum_{j=1}^n \lambda_jx_j^2=C, (x_1, \cdots,
x_n)\in \R^n\}\] for some constant $C$.

  It is not hard to check that $L$ is
 Lagrangian in $\C^n$ with Lagrangian angle given by
 $\theta=(\sum_{j=1}^n \lambda_j)s+\f{\pi}{2}$. It follows that
  $L$
 is  special Lagrangian if $\sum_{j=1}^n\lambda_j=0$.
 In general, a direct computation shows that the induced metric on $L$ is
  independent of $s$.  Hence $\Delta_{L}
 \theta=0$ and $L$ is Hamiltonian stationary.

 Such special Lagrangians were studied by M. Haskins in \cite{ha1}
 \cite{ha2} (for $n=3$) and
 D. Joyce in \cite{jo1} (for general
 dimensions).
 We are informed by D.~Joyce that the Hamiltonian
 stationary ones
 may also be obtained by applying his method
 of ``perpendicular symmetries " in \cite{jo2}.

When $C=0$, the examples are Hamiltonian stationary cones, which generalize the two-dimensional Schoen-Wolfson cones. We
will study the geometry of these examples in the next subsection. Now
assume the constant $C$ in the defining equation is nonzero. 

If $\sum_{j=1}^n\lambda_j\neq 0$, a direct computation shows that 
\[F^\perp=\f{-C}{\sum_{j=1}^n \lambda_j}H.\]  That is,  the submanifold $L$
is a Hamiltonian stationary self-similar solution of the mean
curvature flow.

We summarize the calculations in this subsection in the following
proposition.
\begin{pro}
For any $n$ nonzero integers $\lambda_1,\cdots, \lambda_n$,
 consider the submanifold $L$ of $\C^n$ defined by
\[\{(x_1 e^{i\lambda_1s}, \cdots , x_n e^{i\lambda_ns})\,|
\,\, 0\leq s < 2\pi,\,\sum_{j=1}^n \lambda_jx_j^2=C, (x_1, \cdots,
x_n)\in \R^n\}\] for some constant $C$.  It is special Lagrangian
when $\sum_{j=1}^n\lambda_j=0$. If $\,\,\sum_{j=1}^n\lambda_j\neq
0$, it is Hamiltonian stationary and the normal projection of  the
position vector satisfies
 \[F^\perp=\f{-C}{\sum_{j=1}^n \lambda_j}H.\]  Without
 loss of generality, we can assume that $\sum_{j=1}^{n} \la _{j}>0$
 in this case.
 Then $L$  is a Hamiltonian stationary self-shrinker when $C>0$, a
 Hamiltonian stationary cone when $C=0$, and
 a Hamiltonian stationary self-expander when $C<0$. Moreover, the
Hamiltonian stationary self-expander and the self-shrinker are
asymptotic to the Hamiltonian stationary cone.
\end{pro}
\begin{rem}
The same results hold for any real numbers $\lambda_1,\cdots,
\lambda_n$ which are not all zeros.
\end{rem}
\subsection{The geometry of the examples}
Denote \[F(x_1,\cdots,x_n,s)= (x_1 e^{i\lambda_1s}, \cdots , x_n
e^{i\lambda_ns}).\] It is easy to see that
\[F(x_1,\cdots,x_n,s+\pi)
=F((-1)^{\lambda_1}x_1,\cdots,(-1)^{\lambda_n}x_n,s).\]
As a result, the map $F$ is generically two-to-one.  To solve
this problem, we  restrict the domain to $0 \leq s
<\pi$. The tangent planes at $F(x_1,\cdots,x_n,s+\pi)$ and
$F((-1)^{\lambda_1}x_1,\cdots,(-1)^{\lambda_n}x_n,s)$ agree. They
are both spanned by the vector
\[(i\lambda_1(-1)^{\lambda_1}x_1,\cdots,i\lambda_n(-1)^{\lambda_n}x_n)\]
and the $(n-1)$ plane in $\R^n$ which is perpendicular to the
vector
\[(\lambda_1(-1)^{\lambda_1}x_1,\cdots,\lambda_n(-1)^{\lambda_n}x_n).\]
Define a diffeomorphism $\psi:\R^n \ra \R^n$ by
\[\psi(x_1,\cdots,x_n)=((-1)^{\lambda_1}x_1,\cdots,(-1)^{\lambda_n}x_n).\]
 and a submanifold $\Sigma$ in $\R^n$ by
 \[\Sigma=\{(x_1,\cdots,x_n):
 \sum_{j=1}^{n}\lambda_jx_j^2=C\}.\]
 The restriction of $\psi$ on $\Sigma$ is a diffeomorphism onto
 itself.
 The restriction map is an orientation preserving map if and only if
$\sum_{j=1}^{n}\lambda_{j}$ is even.   It follows that for any $n$
nonzero integers $\lambda_1,\cdots, \lambda_n$, the Lagrangian
submanifold $L'$
\[\{(x_1 e^{i\lambda_1s}, \cdots , x_n e^{i\lambda_ns})\,|
\,\, 0\leq s < \pi,\,\sum_{j=1}^n \lambda_jx_j^2=C, (x_1, \cdots,
x_n)\in \R^n\}\] is oriented if and only if
$\sum_{j=1}^{n}\lambda_{j}$ is even. Assume that $\la_j>0$ for
$1\leq j\leq k$ and   $\la_j<0$ for $k+1\leq j\leq n$. Then the
ansatz of our examples can be written as
\[ \sum_{j=1}^{k}|\la_j|x_j^2=\sum_{j=k+1}^{n}|\la_j|x_j^2+C.\]
The topology of $L'$ is  $ \R^k\times S^{n-k-1}\times S^1$ when
$C<0$, is $S^{k-1}\times \R^{n-k}\times S^1$
 when $C>0$, and is a cone with link $S^{k-1}\times
 S^{n-k-1}\times S^1$, with an isolated singular
point at\/~$0$ when $C=0$.

 To avoid other possible
self-intersections, we require that
$|\lambda_{1}|,\cdots,|\lambda_{n}|$ are pairwise co-prime. This
condition is sufficient to guarantee that the Hamiltonian stationary
cone ($C=0$) is embedded. In addition, we  require that $\la_j=1$ for $1\leq j\leq k$ when $C>0$ and  $\la_j=-1$ for $k+1\leq j\leq n$ in
order to obtain an embedded $L'$ in each case.

Only when $C>0$, $k=1$, or $C<0$, $k=n-1$, the real hypersurface
$\Sigma$ is disconnected. However, if $\la_1$ is odd in the first
case, the factor $e^{i\la _1 s}$ becomes $-1$ at $s=\pi$ and thus
$L'$ will still be connected.  If $\la _n $ is odd in the second
case, the Lagrangian submanifold $L'$ will also be connected for
the same reason.

We summarize these discussions in the following
proposition.
\begin{pro}\label{geoprop}
For any constant $C$ and $n$ nonzero integers $\lambda_1,\cdots,
\lambda_n$,
 the Lagrangian submanifold $L'$ defined by
\[\{(x_1 e^{i\lambda_1s}, \cdots , x_n e^{i\lambda_ns})\,|
\,\, 0\leq s < \pi,\,\sum_{j=1}^n \lambda_jx_j^2=C, (x_1, \cdots,
x_n)\in \R^n\}\] is oriented if and only if
$\sum_{j=1}^{n}\lambda_{j}$ is even.

Assume that $\la_j>0$ for $1\leq j\leq k$ and   $\la_j<0$ for
$k+1\leq j\leq n$. The topology of $L'$ is  $ \R^k\times
S^{n-k-1}\times S^1$ when $C<0$, is $S^{k-1}\times \R^{n-k}\times
S^1$
 when $C>0$, and is a cone with link $S^{k-1}\times
 S^{n-k-1}\times S^1$ when $C=0$.
 If  $C>0$, $k=1$ and $\la _1$ is even, or  $ C<0$, $k=n-1$ and
 $\la _n$ is even,
 there are two connected
 components in $L'$. The submanifold $L'$ is
 connected  for all other cases.

Suppose that $|\lambda_{1}|,\cdots,|\lambda_{n}|$ are pairwise
co-prime. Then the corresponding cones in the case $C=0$ are
embedded. However, one also needs to  require $\la_j=1$ for $1\leq
j\leq k$ to make $L'$ embedded in the case $C>0$, and require
$\la_j=-1$ for $k+1\leq j\leq n$ to make $L'$ embedded in
the case $C<0$.
\end{pro}
\begin{rem}
{\rm Theorem \ref{mainthm}} holds without these extra assumptions
on $\la_j$.
\end{rem}
\begin{rem}
It is worth noting that the case $\sum_{j=1}^{n}\la_j=0$
corresponds to special Lagrangian. Hence the proposition shows
that there are two families of smooth special Lagrangians which
have different topologies but converge to the same special
Lagrangian cone $(C=0)$. The element in one family has topology $
\R^k\times S^{n-k-1}\times S^1$  (for the case $C<0$) and the
element in the other family has topology $S^{k-1}\times
\R^{n-k}\times S^1$ (for the case $C>0$).
\end{rem}
\section{Proof of Theorem \ref{mainthm}}

Assume $\lambda_j$ are non-zero integers and
$\sum_{j=1}^n\lambda_j>0$. Define
\[V_t=\{(x_1 e^{i\lambda_1s}, \cdots , x_n e^{i\lambda_ns})|
\,\, 0\leq s < \pi,\,\sum_{j=1}^n \lambda_jx_j^2=(-2t)\sum_{j=1}^n
\lambda_j, (x_1, \cdots, x_n)\in \R^n \}\] The varifold $V_t$ for
$t \ne 0$ is smooth, and $V_{0}$ is a cone with an isolated
singularity at the origin. As discussed in the previous sections,
$V_t$ are Hamiltonian stationary self-shrinkers for $t<0$ and
Hamiltonian stationary self-expanders for $t>0$. As $t\ra 0$,
$V_t$ converges to the Hamiltonian stationary cone $V_0$. The
geometry of $V_t$ is discussed in Proposition \ref{geoprop}. What
is left in the proof of Theorem \ref{mainthm} is to show that the
varifolds $V_t$ for $-\infty<t<\infty$ form an eternal solution of
Brakke flow without mass loss. We first recall the definition of
Brakke flow.

\subsection{Brakke flow}\label{brakke}
 A family of varifolds
$V_{t}$ is said to form a solution of the Brakke flow \cite{br} if
\begin{equation}\label{brakkeineq}\bar{D}|\!|V_{t}|\!|(\phi)\leq
\delta (V_{t},\phi)(h(V_{t}))\end{equation}
 for each $\phi \in
C^{1}_{0}(\R^{n})$ with $\phi \geq 0$, where
$\bar{D}|\!|V_{t}|\!|(\phi)$ is the upper derivative defined by $
\overline{\lim}_{t_{1}\rightarrow t }
\f{|\!|V_{t_{1}}|\!|(\phi)-|\!|V_{t}|\!|(\phi)} {t_{1}-t}$ and
$h(V_{t})$ is the generalized mean curvature vector of $V_{t}$. In
the setting of this paper,
$$ \delta (V_t,\phi)(h(V_t))=-\int \phi |h(V_t)|^{2} d|\!|V_t|\!|
+\int D\phi\cdot h(V_t)d|\!|V_t|\!|.$$

In our case, the family $V_t$ satisfy mean curvature flow for
$t<0$ and $t>0$ and
 the singularity only happens at the $t=0$ slice. The
following proposition is formulated in \cite{lw} as a criterion to
check the solutions of Brakke flow in this situation.

\begin{pro}\label{mainpro}

Suppose the varifold $V_t$, $a <t<b$ forms a smooth mean curvature
flow in $\R^n$ except at $t=c\in (a, b)$ and $|\!|V_t|\!|$
converges in Radon measure to $|\!|V_c|\!|$ as $t\rightarrow c$.
If $ \lim_{t\rightarrow c^-} \f{d}{dt}|\!|V_t|\!|(\phi)$ and
$\lim_{t\rightarrow c^+} \f{d}{dt}|\!|V_t|\!|(\phi)$ are both
either finite or $-\infty$ and

\begin{equation}\label{assumption}\lim_{t\rightarrow
c^\pm} \f{d}{dt}|\!|V_t|\!|(\phi)\leq \delta (V_0,
\phi)(h(V_0))\end{equation} for any $\phi\in C^1_0(\R^n)$ then
$V_t$ forms a solution of the Brakke flow.
\end{pro}
\begin{dfn}
If   $V_t$ form a solution of the Brakke flow for
$-\infty<t<\infty$, we call it an eternal solution for Brakke
flow. Moreover, if the equality  in {\rm (\ref{brakkeineq})} is
achieved for all $-\infty<t<\infty$, we say the solution has no
mass loss.
\end{dfn}
\subsection{Completion of the proof}
Since for a smooth mean curvature flow, we have

\[\f{d}{dt}|\!|V_t|\!|(\phi)=\delta (V_t,\phi)(h(V_t))=-\int \phi
|h(V_t)|^{2} \dvt +\int D\phi \cdot h(V_t) \dvt.\] To
 apply Proposition \ref{mainpro} and prove that
the equality in the Brakke flow is achieved, it suffices to show

\begin{equation}\label{mainid-}\begin{split}
&\lim_{t\rightarrow 0^-}-\int \phi \,|h(V_t)|^{2}
\dvt
+\int D\phi \cdot h(V_t) \dvt\\
&=-\int \phi \,|h(V_0)|^{2} \dvzero +\int D\phi \cdot h(V_0)
\dvzero,\end{split}\end{equation}
and
\begin{equation}\label{mainid+}\begin{split}
&\lim_{t\rightarrow 0^+}-\int \phi \,|h(V_t)|^{2}
\dvt
+\int D\phi \cdot h(V_t) \dvt\\
&=-\int \phi \,|h(V_0)|^{2} \dvzero +\int D\phi \cdot h(V_0)
\dvzero.\end{split}\end{equation}
 A direct calculation shows that \begin{equation}\label{etnorm2}
 \left|V_t\right|^2= \sum_{j=1}^{n}
x_j^2, \end{equation}
\begin{equation}\label{hetnormsq2}\left|h(V_t)\right|^2
=\frac{\sum_{j=1}^{n}
\lambda_j^2}{\sum_{j=1}^{n}\lambda_j^2x_j^2},\end{equation} and

\begin{equation}\label{dett2}d|\!|V_t|\!|=
\sqrt{\sum_{j=1}^{n}\lambda_j^2x_j^2} \,\,dS_t \,ds
,\end{equation} where $dS_t$ is the volume form of the
hypersurface
\[\Sigma_t= \{(x_1, \cdots , x_n )|
\sum_{j=1}^n \lambda_jx_j^2=(-2t)\sum_{j=1}^n \lambda_j \}\] in $\R^n$.

We can parameterized $\Sigma_t$ by rewriting the defining equation
as

\[\sum_{j=1}^k|\lambda_j|x_j^2=\sum_{j=k+1}^n|\lambda_j| x_j^2
-2t\sum_{j=1}^n\lambda_j\] where $\lambda_j>0$ for $j=1\cdots k$
and $\lambda_j<0$ for $j=k+1,\cdots, n$.

Suppose $X_2=(0,\cdots, 0, x_{k+1},\cdots, x_n)$ gives the
embedding of the surface $\sum_{j=k+1}^n|\lambda_j| x_j^2=1$ and
$X_1=(x_1, \cdots, x_k,0, \cdots, 0)$ gives the embedding of the
surface $\sum_{j=1}^k|\lambda_j|x_j^2=1$. Then the hypersurface
$\Sigma_t$ for $t<0$ can be parameterized by

\[X=\left(r^2-2t\sum_{j=1}^n\lambda_j\right)^{\frac{1}{2}}X_1+rX_2,\]
where
 $r^2=\sum_{j=k+1}^n|\lambda_j| x_j^2$.

It is not hard to check that the volume form of $\Sigma_t$, $t<0$
is given by

\begin{equation}\label{volform1}dS_t=r^{n-k-1} (r^2-2t\sum_{j=1}^n \lambda_j)^{\frac{k-1}{2}}
\left(\frac{r^2}{r^2-2t\sum_{j=1}^n
\lambda_j}|X_1^\perp|^2+|X_2^\perp|^2\right)^{\frac{1}{2}}drdS_t^-
dS_t^+,\end{equation} where
 $dS_t^-$ is the volume form of $$\{(x_{k+1},\cdots,
x_n)|\sum_{j=k+1}^n |\lambda_j|x_j^2=1\}\subset \R^{n-k},$$ and
$dS_t^+$ is the volume form of $\{(x_{1},\cdots, x_k)|\sum_{j=1}^k
|\lambda_j|x_j^2=1\}\subset \R^k$.

From (\ref{etnorm2}), (\ref{hetnormsq2}),
 (\ref{dett2}) and (\ref{volform1}), we have
\begin{equation}\label{mainid-1}\begin{split}&\int \phi \,|h(V_t)|^{2} \dvt\\
= &\int \f{\phi \sum _{j=1}^{n}\lambda_j^2
}{\sqrt{\sum_{j=1}^n\lambda_j^{2}x_j^2}}r^{n-k-1}
 (r^2-2t\sum_{j=1}^n \lambda_j)^{\frac{k-1}{2}}
\left(\frac{r^2}{r^2-2t\sum_{j=1}^n
\lambda_j}|X_1^\perp|^2+|X_2^\perp|^2\right)^{\frac{1}{2}}
drdS_t^- dS_t^+ ds. \end{split}\end{equation}
Because  $t<0$, and
 $\sum_{j=1}^{n}\lambda_j>0 $, it follows that
\begin{equation}\label{volestimate}\begin{split}r^{n-k-1}
 (r^2-2t\sum_{j=1}^n \lambda_j)^{\frac{k-1}{2}}&
 < (2r^2-2t\sum_{j=1}^{n}\lambda_j)^{\frac{n-2}{2}}\\
 &=(\sum_{j=1}^n|\lambda_j|x_j^2)^{\frac{n-2}{2}}\\
 &\leq(\sum_{j=1}^n\lambda_j^{2}x_j^2)^{\frac{n-2}{2}},
 \end{split}\end{equation}
and
 the integrand in (\ref{mainid-1}) is bounded by the function
 \[ \phi \sum_{j=1}^n\la _{j}^{2}\,(\sum_{j=1}^n\lambda_j^{2}x_j^2)^{\frac{n-3}{2}}
 \left(|X_1^\perp|^2+|X_2^\perp|^2\right)^{\frac{1}{2}}.\]
Moreover, the function $\phi$ has compact support and
$\{\la_j\}_{j=1}^{n}$ is fixed, so this is an integrable function
when $n\ge 3$. By the dominate convergence theorem, we thus have
\[\begin{split}
\lim_{t\rightarrow 0^-}-\int \phi \,|h(V_t)|^{2} \dvt &=-\int \phi
\,|h(V_0)|^{2} \dvzero \end{split}.\] The same estimates also show
\[\begin{split} \lim_{t\rightarrow 0^-}\int D\phi \cdot h(V_t)
\dvt =\int D\phi \cdot h(V_0) \dvzero.\end{split}\] We thus prove
(\ref{mainid-}). Note that when $k=1$ or $n-1$,  the expression
above is slightly different, but the same argument works.

When $t>0$, we rewrite the defining equation as

\[\sum_{j=1}^k|\lambda_j|x_j^2+2t\sum_{j=1}^n\lambda_j
=\sum_{j=k+1}^n|\lambda_j| x_j^2 ,\] where $\lambda_j>0$ for
$j=1\cdots k$ and $\lambda_j<0$ for $j=k+1,\cdots, n$.  Then the
hypersurface $\Sigma_t$ for $t>0$ can be parameterized by

\[X=rX_1+\left(r^2+2t\sum_{j=1}^n\lambda_j\right)^{\frac{1}{2}}X_2,\]
where
 $r^2=\sum_{j=1}^k|\lambda_j| x_j^2$.

  Similar computations as in the case $t<0$ show that the volume
  form  of $\Sigma_t$, $t>0$ is given by

\begin{equation}\label{volform2}dS_t=r^{k-1} (r^2+2t\sum_{j=1}^n \lambda_j)^{\frac{n-k-1}{2}}
\left(|X_1^\perp|^2+\frac{r^2}{r^2+2t\sum_{j=1}^n
\lambda_j}|X_2^\perp|^2\right)^{\frac{1}{2}}drdS_t^-
dS_t^+,\end{equation} Therefore,

\begin{equation}\label{mainid+1}\begin{split}&\int \phi \,|h(V_t)|^{2} \dvt\\
= &\int \f{\phi \sum_{j=1}^{n} \lambda_j^2
}{\sqrt{\sum_{j=1}^n\lambda_j^{2}x_j^2}} r^{k-1}
(r^2+2t\sum_{j=1}^n \lambda_j)^{\frac{n-k-1}{2}}
\left(|X_1^\perp|^2+\frac{r^2}{r^2+2t\sum_{j=1}^n
\lambda_j}|X_2^\perp|^2\right)^{\frac{1}{2}} drdS_t^- dS_t^+ ds.
\end{split}\end{equation}

Because  $t>0$, and
 $\sum_{j=1}^{n}\lambda_j>0 $, we can similarly show that
the integrand in (\ref{mainid+1}) is bounded by the function
 \[ \phi \sum_{j=1}^{n}\la _{j}^{2}\,(\sum_{j=1}^n\lambda_j^{2}x_j^2)^{\frac{n-3}{2}}
 \left(|X_1^\perp|^2+|X_2^\perp|^2\right)^{\frac{1}{2}},\]
which is an integrable function if  $n\ge 3$ and $\phi$ has
compact support. By the dominate convergence theorem, we thus have
\[\begin{split}
\lim_{t\rightarrow 0^+}-\int \phi \,|h(V_t)|^{2} \dvt &=-\int \phi
\,|h(V_0)|^{2} \dvzero \end{split}.\] The same estimates also show
\[\begin{split} \lim_{t\rightarrow 0^+}\int D\phi \cdot h(V_t)
\dvt =\int D\phi \cdot h(V_0) \dvzero.\end{split}\] We thus prove
(\ref{mainid+}). Again, when $k=1$ or $n-1$, the expression above
needs slight modification, but the same argument gives the
conclusion. Since the two dimensional case is already proved in
\cite{lw}, this completes the proof of Theorem \ref{mainthm}.

\section{Proof of Theorem \ref{thm2}}
Recall \cite[Theorem A]{jlt} that if
\begin{equation}
\begin{aligned}
\frac{d w_j}{d s}&= \la_je^{i\th(s)}\,\ov{w_1\cdots w_{j-1}
w_{j+1}\cdots w_n},\qquad j=1,\cdots,n,\\
\frac{d\th}{d s}&=\al\Im(e^{-i \th(s)} w_1\cdots w_n),
\end{aligned}
\label{lm4eq1}
\end{equation}
then the submanifold\/ $L$ in $\C^n$ given by
\begin{equation*}
L=\bigl\{\bigl(x_1w_1(s),\cdots,x_nw_n(s)\bigr):\text{$s\in I,$
$x_1,\cdots,x_n\in\R,$ $\ts\sum_{j=1}^{n}\la_jx_j^2=C$}\bigl\},
\end{equation*}
is Lagrangian, with Lagrangian angle $\th(s)$ at\/
$(x_1w_1(s),\cdots,x_nw_n(s)),$ and its position vector\/ $F$ and
mean curvature vector\/ $H$ satisfy\/ $\al F^\perp=CH$.

Because there  is a lot of freedom to rescale the constants (see
\cite[Remark 3.2]{jlt}), we can assume $\al=1$ for simplicity.
From Theorem~F of \cite{jlt}, there is a dense set of initial data
such that the solutions $w_1(s),\cdots,w_n(s)$ and $\th(s)$ of
\eq{lm4eq1} are periodic. Suppose the period is $T$. Then  $V_t$,
which is defined by
\[ V_t=\{(x_1 w_1(s), \cdots ,
x_n w_n(s)):\, 0\leq s < T,\, \sum_{j=1}^n \lambda_jx_j^2=2t,
(x_1, \cdots, x_n)\in \R^n \},\] is  a Lagrangian self-shrinker
for $t<0$, a Lagrangian self-expander for $t>0$, and a Lagrangian
cone for $t=0$.

To show this family form an eternal solution of Brakke flow
without mass loss, like the proof of Theorem~\ref{mainthm}, it
remains to check \eq{mainid-} and \eq{mainid+}.

Denote $w _j(s)=r_{j}(s)e^{i\,\varphi _j(s)}$ and
$\varphi(s)=\sum_{j=1}^{n} \varphi _j(s)$, where $r_j(s)=|w
_j(s)|$. Remember that $r_j$ has positive lower and upper
bounds.
 A direct calculation shows that \begin{equation}\label{getnorm2}
 \left|V_t\right|^2= \sum_{j=1}^{n}r_j^2
x_j^2, \end{equation}
\begin{equation}\label{ghetnormsq2}\left|h(V_t)\right|^2
=\bigl(\sum_{j=1}^{n}\frac{\lambda_j^2x_j^2}{r_j^2}\bigr)^{-1}\sin^{2}
(\varphi-\th ),\end{equation} and

\begin{equation}\label{gdett2}d|\!|V_t|\!|=\frac{r_1^2\cdots r_n^2}{\sqrt{\sum_{j=1}^{n}\la_j^2
x_j^2}}\sum_{j=1}^{n}\frac{\la_j^2x_j^2}{r_j^2}\,\,dS_t \,ds
,\end{equation} where $dS_t$ is the volume form of the
hypersurface
\[\Sigma_t= \{(x_1, \cdots , x_n )|
\sum_{j=1}^n \lambda_jx_j^2=2\,t \}\] in $\R^n$. When $t<0$, we
can parameterized $\Sigma_t$ by rewriting the defining equation as

\[\sum_{j=1}^k|\lambda_j|x_j^2-2\,t=\sum_{j=k+1}^n|\lambda_j| x_j^2
.\]
Suppose $X_1=(x_1, \cdots, x_k,0, \cdots, 0)$ gives the
embedding of the surface $\sum_{j=1}^k|\lambda_j|x_j^2=1$ and
$X_2=(0,\cdots, 0, x_{k+1},\cdots, x_n)$ gives the embedding of
the surface $\sum_{j=k+1}^n|\lambda_j| x_j^2=1$. Then the
hypersurface $\Sigma_t$ for $t<0$ can be parameterized by

\[X=rX_1+(r^2-2\,t)^{\frac{1}{2}}X_2,\]
where
 $r^2=\sum_{j=1}^k|\lambda_j| x_j^2$.

It is not hard to check that the volume form of $\Sigma_t$, $t<0$
is given by
\begin{equation}\label{gvolform1}dS_t=r^{k-1} (r^2-2\,t)^{\frac{n-k-1}{2}}
\left(|X_1^\perp|^2+\frac{r^2}{r^2-2\,t}|X_2^\perp|^2\right)^{\frac{1}{2}}drdS_t^-
dS_t^+,\end{equation} where
 $dS_t^-$ is the volume form of $$\{(x_{k+1},\cdots,
x_n)|\sum_{j=k+1}^n |\lambda_j|x_j^2=1\}\subset \R^{n-k},$$ and
$dS_t^+$ is the volume form of $\{(x_{1},\cdots, x_k)|\sum_{j=1}^k
|\lambda_j|x_j^2=1\}\subset \R^k$.

From (\ref{getnorm2}), (\ref{ghetnormsq2}),
 (\ref{gdett2}) and (\ref{gvolform1}), we have
\begin{equation}\label{gmainid-1}\begin{split}&\int \phi \,|h(V_t)|^{2} \dvt\\
= &\int \frac{\phi \sin^2(\varphi-\th)r_1^2\cdots
r_n^2}{\sqrt{\sum_{j=1}^{n}\la_j^2 x_j^2}}r^{k-1}
 (r^2-2t)^{\frac{n-k-1}{2}}
\left(|X_1^\perp|^2+\frac{r^2}{r^2-2t}|X_2^\perp|^2\right)^{\frac{1}{2}}
drdS_t^- dS_t^+ ds.
\end{split}\end{equation}
Similar to \eq{volestimate}, we have
\begin{equation}\label{gvolestimate}r^{k-1}
 (r^2-2t)^{\frac{n-k-1}{2}}<
 \bigl(\sum_{j=1}^n|\lambda_j|x_j^2\bigr)^{\frac{n-2}{2}}
 \leq\bigl(\sum_{j=1}^n\lambda_j^{2}x_j^2\bigr)^{\frac{n-2}{2}}.
 \end{equation}
Here for simplicity we use the normalization in \cite[Remark
3.2]{jlt} and assume $|\la_j|\ge 1$ for all $j$.
 Because $\phi$ is a $C^1$
function with compact support, $t<0$, and $r_1,\cdots,r_{n}$ are
bounded, the integrand in (\ref{gmainid-1}) is bounded by the
function $$C(\sum_{j=1}^n\lambda_j^{2}x_j^2)^{\frac{n-3}{2}}
 \left(|X_1^\perp|^2+|X_2^\perp|^2\right)^{\frac{1}{2}}.$$
 This is an integrable function
when $n\ge 3$. By the dominate convergence theorem, we thus have
\[\begin{split}
\lim_{t\rightarrow 0^-}-\int \phi \,|h(V_t)|^{2} \dvt &=-\int \phi
\,|h(V_0)|^{2} \dvzero \end{split}.\] The same estimates also show
\[\begin{split} \lim_{t\rightarrow 0^-}\int D\phi \cdot h(V_t)
\dvt =\int D\phi \cdot h(V_0) \dvzero.\end{split}\] We thus prove
(\ref{mainid-}) when $n\ge 3$.  Again, when $k=1$ or $n-1$, the
expression above needs slight modification, but the same argument works.

When $n=2$, a direct computation gives
\begin{equation}\label{2did1}\begin{split}\int \phi \,|h(V_t)|^{2} \dvt&=
\int \frac{\phi \sin^2(\varphi-\th)r_1^2
r_2^2}{\sqrt{\sum_{j=1}^{2}\la_j^2 x_j^2}}
\left(1+\frac{x_1^2}{|\la_2|\,x_2^2}\right)^{\frac{1}{2}}dx_1 ds,
\end{split}\end{equation}
where $|\la_1|x_1^2-2\,t=|\la_2|x_2^2$. We proceed as in \cite{lw}
by dividing into two cases $\phi(0)=0$ or $\phi(0)\ne 0$. When
$\phi (0)=0$, we have $\phi (V_t)\le C|V_t|\le
C'(\sum_{j=1}^2\lambda_j^{2}x_j^2)^{\frac{1}{2}}.$ With this extra
power, the integrand becomes bounded  and (\ref{mainid-}) follows
from the dominate convergence theorem.

 If
$\phi(0)\ne 0$, we have $\int \phi \,|h(V_0)|^{2} \dvzero=\infty$.
We can also prove that  $$\lim_{t\rightarrow 0^-}\int \phi
\,|h(V_t)|^{2} \dvt=\infty$$ and both $\lim_{t\rightarrow 0^-}\int
D\phi \cdot h(V_t) \dvt$ and $\int D\phi \cdot h(V_0) \dvzero$ are
finite. Hence (\ref{mainid-}) holds trivially. Now we give the
proof for these facts. Given any $\ep>0$, there exist a $\delta>0$
such that $|\phi|\ge \ep$ in $B_\delta(0)$. Using the
normalization in \cite[Remark 3.2]{jlt}, we can  assume $\la_1=1$
and $\la_2=-1$ for simplicity. Thus $x_1^2-2t=x_2^2$ and because
$r_1$, $r_2$ and $\sin^2(\varphi-\th)$ all have positive lower
bound, from \eq{2did1} and $0\le s\le T$ we have
\begin{equation}\label{2did2}\begin{split}\int \phi \,|h(V_t)|^{2}
\dvt&\ge C\ep\! \int_0^a \frac{1}{\sqrt{2 x_1^2-2t}} dx_1.
\end{split}\end{equation}
Note that  $|V_t|^2=x_1^2r_1^2+x_2^2r_2^2\le \de^2$ implies
$x_1^2(r_1^2+r_2^2)-2tr_2^2\le \de^2$. Hence  when $t$ is close to
$0^-$, this set  contains a uniform interval $[0,a]$ where $a$
depends only on $\de$. Since $\sqrt{2x_1^2-2t}\le
\sqrt{2}|x_1|+\sqrt{2}\sqrt{-t}$,  \eq{2did2} becomes
\begin{equation}\label{2did3}\begin{split}\int \phi \,|h(V_t)|^{2}
\dvt&\ge C'\ep\ln \frac{a+\sqrt{-t}}{\sqrt{-t}},
\end{split}\end{equation} which tends to $\infty$ as $t$ tends to $0^-$.
A direct computation gives $$\int D\phi \cdot h(V_t) \dvt \le \int
\frac{|D\phi|\, r_1^2 r_2^2}{\sqrt{2 x_1^2-2t}}\sqrt{
\frac{x_1^2}{r_l^2}+\frac{x_1^2-2t}{r_2^2}}
\left(1+\frac{x_1^2}{x_1^2-2t}\right)^{\frac{1}{2}}dx_1 ds.$$ The
integrand is bounded and hence we can use the dominate convergence
theorem to show that the limit is finite. This prove \eq{mainid-}
when $n=2$

When $t>0$, similar arguments give \eq{mainid+} for $n\ge3$ and
$n=2$. Thus Theorem \ref{thm2} is proved.
\begin{rem}
Similar to the discussions in {\rm \S2.2 (}and also see {\rm
\cite[Theorem~A and Theorem~F]{jlt})}, $V_t$ is a closed,
nonsingular, immersed Lagrangian self-expander in $\C^n$
diffeomorphic to $S^{k-1}\times \R^{n-k}\times S^1$ when $t>0,$
and a closed, nonsingular, immersed Lagrangian self-shrinker in
$\C^n$ diffeomorphic to $\R^{k} \times S^{n-k-1}\times S^1$ when
$t<0,$ and\/ $V_0$ is a closed, immersed Lagrangian cone in $\C^n$
with link\/ $S^{k-1}\times S^{n-k-1}\times S^1,$ with an isolated
singular point at\/~$0$. To study the embeddedness of these examples,
one needs to have a better understanding on $\ga_j$, which is
defined in {\rm \cite[Theorem~E]{jlt}}.
\end{rem}

\end{document}